\documentclass{ifacconf}
\usepackage[round]{natbib} 
\usepackage{graphicx}      

\usepackage{amsmath,amssymb}
\usepackage{url}
\usepackage{diffcoeff}
\usepackage{bm}

\usepackage{mathrsfs}
\usepackage{multirow}

\usepackage{caption, subfig}

\usepackage{xcolor,colortbl}

\newtheorem{remark}{Remark}
\newtheorem{conjecture}{Conjecture}

\DeclareMathOperator*{\esssup}{ess\,sup}
\DeclareMathOperator{\Tr}{Tr}

\def\onedot{$\mathsurround0pt\ldotp$}
\def\cddot{
	\mathbin{\vcenter{\baselineskip.67ex
			\hbox{\onedot}\hbox{\onedot}}%
}}

\makeatletter \renewcommand\d[1]{\ensuremath{%
		\;\mathrm{d}#1\@ifnextchar\d{\!}{}}}
\makeatother

\graphicspath{{./Figures/}}
\begin{document}
\begin{frontmatter}

\title{Structure-preserving discretization of port-Hamiltonian plate models \thanksref{footnoteinfo}} 

\thanks[footnoteinfo]{This work is  supported by the project ANR-16-CE92-0028,
	entitled {\em Interconnected Infinite-Dimensional systems for Heterogeneous
		Media}, INFIDHEM, financed by the French National
	Research Agency (ANR) and the Deutsche Forschungsgemeinschaft (DFG). Further information is available at {\url{https://websites.isae-supaero.fr/infidhem/the-project}}.
	}

\author[ISAE]{Andrea Brugnoli}
\author[ISAE]{Daniel Alazard} 
\author[ISAE]{Val\'erie Pommier-Budinger}
\author[ISAE]{Denis Matignon}

\address[ISAE]{ISAE-SUPAERO, Universit\'e de Toulouse, France.\\
	10 Avenue Edouard Belin, BP-54032, 31055 Toulouse Cedex 4. \\
	Andrea.Brugnoli@isae.fr,  Daniel.Alazard@isae.fr, \\
	Valerie.Budinger@isae.fr, Denis.Matignon@isae.fr}

\begin{abstract}
Methods for discretizing port-Hamiltonian systems are of interest both for simulation and control purposes. Despite the large literature on mixed finite elements, no rigorous analysis of the connections between mixed elements and port-Hamiltonian systems has been carried out. In this paper we demonstrate how existing methods can be employed to discretize dynamical plate problems in a structure-preserving way. Based on convergence results of existing schemes, new error estimates are conjectured; numerical simulations confirm the expected behaviors.
\end{abstract}


\end{frontmatter}

\section{Introduction}
Distributed port-Hamiltonian (dpH) systems, since their introduction in \cite{VANDERSCHAFTdph}, have attracted a lot of attention. For simulation and control design purposes, a suitable, i.e. structure preserving, discretized model is needed. When dealing with higher geometrical dimensions, obtaining a finite-dimensional approximation is not an easy task. One possible strategy is to make use of finite element discretization. In \cite{WeakForm_Kot} the authors provide a way to discretize systems of conservation laws using finite element exterior calculus. However, this methodology does not easily generalize to more complicated models (e.g. elasticity problems).

Thanks to \cite{CardosoRibeiro2018}, it has become evident that there is a strict link between  discretization of port-Hamiltonian (pH) systems and mixed finite elements. Velocity-stress formulation for the wave dynamics and elastodynamics problems are indeed Hamiltonian and their mixed discretization preserves such a structure (cf. \cite{Kirby2015}, where a symplectic in time and space scheme is constructed for the wave equation). This allows using known finite element scheme to preserve the pH structure at the discrete level.

Mixed finite elements for the wave equation have been studied in \cite{Geveci,becacheWave}. For elastodynamics the construction of stable elements gets more complicated because of the presence of the symmetric stress tensor. Existing elements enforce symmetry either strongly (\cite{becacheElas}) or weakly (\cite{ArnoldWeak,mindlinVeiga}).

In this paper, the mixed finite element discretization of plate models in pH form is studied. These models have been recently presented in \cite{BRUGNOLI2019940,BRUGNOLI2019961}, but without any convergence analysis. Here, available mixed finite elements for the wave dynamics and elastodynamics are adapted to the Mindlin plate problem. Error bounds are conjectured. For the Kirchhoff plate, the differential operator is of second order. Hence, more regular elements are required. The Hellan-Herrmann-Jhonson scheme is used (\cite{Blum1990,arnold2019hellan}). It is conjectured that the convergence results valid for the static problem carry over to the dynamical case. Numerical simulations are implemented to confirm our conjectures.

The paper is organized as follows. In Section \ref{sec:pHplates}, plate models as port-Hamiltonian systems are briefly recalled. In Sec.~\ref{sec:mixed}, the weak formulation and corresponding finite element combinations are illustrated. The discretization relies on existing finite elements, hence the numerical implementation requires little effort. In Sec. \ref{sec:numerics}, the numerical results, which confirm the expected behavior, are presented. The implementation is performed using the Firedrake python library (\cite{rathgeber2017firedrake}).

\section{Plate models in port-Hamiltonian form}
\label{sec:pHplates}

In this section the models under consideration are recalled. More details can be found in  \cite{BRUGNOLI2019940} for the Mindlin plate and in \cite{BRUGNOLI2019961} for the Kirchhoff plate. 

\subsection{Notations}
The space of all, symmetric and skew-symmetric $d\times d$ matrices are denoted by $\mathbb{M}, \mathbb{S}, \mathbb{K}$ respectively. The space of $\mathbb{R}^d$ vectors is denoted by $\mathbb{V}$. $\Omega \subset \mathbb{R}^d$ is an open connected set. The geometric dimension of interest in this paper is $d=2$. For a scalar field $u: \Omega \rightarrow \mathbb{R}$ the gradient is defined as 
\begin{equation*}
\mathrm{grad}(u) =  \nabla u := \begin{pmatrix}
\partial_{x_1} u \dots \partial_{x_d} u \\
\end{pmatrix}^\top.
\end{equation*}
For a vector field $\bm{u}: \Omega \rightarrow \mathbb{V}$, with components $u_j$, the gradient (Jacobian) is defined as
\begin{equation*}
\mathrm{grad}(\bm{u})_{i j}:= (\nabla \bm{u})_{ij} = \partial_{x_j} u_i.
\end{equation*}
The symmetric part of the gradient operator $\mathrm{Grad}$ (i. e. the deformation tensor in continuum mechanics) is thus given by
\begin{equation*}
\mathrm{Grad}(\bm{u}) := \frac{1}{2} \left(\nabla \bm{u} + \nabla^\top \bm{u} \right).
\end{equation*}
The Hessian operator of $u$ is then computed as follows
\begin{equation*}
\mathrm{Hess}(u) = \nabla^2 u = \mathrm{Grad}(\mathrm{grad}(u)),
\end{equation*}
For a tensor field $\bm{U}: \Omega \rightarrow \mathbb{M}$, with components $u_{ij}$, the divergence is a vector, defined column-wise as
\begin{equation*}
\mathrm{Div}(\bm U) = \nabla \cdot \bm{U} := \left( \sum_{i = 1}^d \partial_{x_i} u_{ij} \right)_{j = 1, \dots, d}.
\end{equation*}
The double divergence of a tensor field $\bm{U}$ is then a scalar field defined as
\begin{equation*}
\mathrm{div}(\mathrm{Div}(\bm U)):= \sum_{i, j = 1}^d \partial_{x_i} \partial_{x_j} u_{ij}.
\end{equation*}
The $L^2$ inner products of scalar, vector and matrix fields are defined as
\begin{align*}
	(u, v) &= \int_{\Omega} u \ v \d\Omega, \quad u, v : \Omega \rightarrow \mathbb{R}, \\
	(\bm{u}, \bm{v}) &= \int_{\Omega} \bm{u} \cdot \bm{v} \d\Omega, \quad \bm{u}, \bm{v} : \Omega \rightarrow \mathbb{V}, \\
	(\bm{U}, \bm{V}) &= \int_{\Omega} \bm{U} \cddot \bm{V} \d\Omega, \quad \bm{U}, \bm{V} : \Omega \rightarrow \mathbb{M},
\end{align*}
where $\bm{u} \cdot \bm{v} := \sum_{i} u_{i} v_{i}$ is the scalar product in $\mathbb{V}$ and $\bm{U} \cddot \bm{V} := \sum_{i,j} u_{ij} v_{ij}$ is the tensor contraction. The standard notation $H^m(\Omega)$ denotes the Sobolev space of square integrable functions with  $m^\text{th}$ derivative in $L^2$ and norm $||\cdot||_m$. In particular, $H^1_0(\Omega)$ is the space of weakly derivable functions with vanishing trace. For $\mathbb{X} \subseteq \mathbb{M}$, let
\begin{equation*}
\begin{aligned}
H(\mathrm{div}, \Omega) &= \{\bm{u} \in L^2(\Omega, \mathbb{V}) \vert \; \mathrm{div}(\bm{u}) \in L^2(\Omega) \}, \\
H(\mathrm{Div}, \Omega; \mathbb{X}) &= \{\bm{U} \in L^2(\Omega, \mathbb{X}) \vert \; \mathrm{Div}(\bm{U}) \in L^2(\Omega; \mathbb{V}) \},
\end{aligned}
\end{equation*}
which are Hilbert spaces with the norm $||\bm{u}||^2_{\text{div}} = ||\bm{u}||^2 + ||\mathrm{div}(\bm{u})||^2, \; ||\bm{U}||^2_{\text{Div}} = ||\bm{U}||^2 + ||\mathrm{Div}(\bm{U})||^2$. The following abbreviations will be used
\begin{equation*}
\begin{aligned}
M &= H(\mathrm{Div}, \Omega; \mathbb{M}), \\
S &= H(\mathrm{Div}, \Omega; \mathbb{S}),
\end{aligned} \qquad
\begin{aligned}
D &= H(\mathrm{div}, \Omega), \\
L &= L^2(\Omega),
\end{aligned} \qquad
\begin{aligned}
V &= L^2(\Omega; \mathbb{V}), \\
K &= L^2(\Omega; \mathbb{K}).
\end{aligned}
\end{equation*}
Let $\mathcal{X}$ be a Hilbert space, and $t_f$ a positive real number. We denote by $L^\infty([0, t_f]; \mathcal{X})$ or $L^\infty(\mathcal{X})$ the space of functions $f: [0, t_f] \rightarrow X$ for which the time-space norm $||\cdot||_{L^\infty([0, t_f]; \mathcal{X})}$ satisfies
\[
||f||_{L^\infty([0, t_f]; \mathcal{X})} = \esssup_{t \in [0,t_f]} ||f||_{\mathcal{X}} < \infty.
\]

\subsection{Mindlin-Reissner plate}
The Mindlin model is a generalization to the 2D case of the Timoshenko beam model and is expressed by a system of two coupled PDEs (\cite{timoshenko1959theory}) 
\begin{equation}
\label{eq:clMin}
\begin{cases}
\displaystyle \rho b \diffp[2]{w}{t} &= \mathrm{div}(\bm{q}) + f, \quad (\bm{x}, t) \in \Omega \times [0, t_f],  \vspace{1mm}\\
\displaystyle \frac{\rho b^3}{12} \diffp[2]{\bm \theta}{t} &= \bm{q} + \mathrm{Div}(\bm M) + \bm{\tau}, \\
\end{cases}
\end{equation}
where $\rho$ is the mass density, $b$ the plate thickness, $w$ the vertical displacement, $\bm \theta = (\theta_x, \theta_y)^\top$ collects the deflection of the cross section along axes $x$ and $y$ respectively. The fields $f, \bm{\tau}$ represent distributed forces and torques. Variables $\bm{M}, \bm{q}$ represent the momenta tensor and the shear stress. Hooke's law relates those to the curvature tensor and shear deformation vector
\begin{equation*}
\begin{aligned}
\bm{M} &:= \mathcal{D} \bm{K} \in \mathbb{S}, \\ \bm{q} &:= \mathcal{C} \bm{\gamma},
\end{aligned} \qquad
\begin{aligned}
\bm{K} &:= \mathrm{Grad}(\bm{\theta}) \in \mathbb{S}, \\ \bm{\gamma} &:= \mathrm{grad}(w) - \bm{\theta},.
\end{aligned}
\end{equation*}
Tensors $\mathcal{D}, \ \mathcal{C}$ are symmetric positive 
\begin{equation}
\label{eq:bend_rig_tensor}
	\mathcal{D} (\cdot) = \frac{E_Y b^3}{12 (1 - \nu^2)}[(1-\nu)(\cdot) + \nu \Tr(\cdot)], \quad \mathcal{C} (\cdot) = \frac{E b k }{2(1+\nu)}(\cdot),
\end{equation}
where $E_Y$ is the Young modulus, $\nu$ is the Poisson modulus, $k$ is the shear correction factor.
 The kinetic and potential energies  $E_c, E_p$ read
\begin{equation}
\begin{aligned}
E_c &=  \frac{1}{2} \int_{\Omega} \left\{ \rho b \left(\diffp{w}{t} \right)^2 +  \frac{\rho b^3}{12} \diffp{\bm{\theta}}{t} \cdot \diffp{\bm{\theta}}{t}  \right\} \d\Omega, \\
E_p &= \frac{1}{2} \int_{\Omega} \left\{ \bm{M} \cddot \bm{K} + \bm{q} \cdot \bm{\gamma}  \right\}\d\Omega.
\end{aligned}
\end{equation} 
The Hamiltonian  is easily written as $H = E_c + E_p$. To get a port-Hamiltonian formulation, suitable energy variables must be selected. The appropriate set is the following:
\begin{equation}
\begin{aligned}
\alpha_w &= \rho b \diffp{w}{t}, \\
\bm{A}_{\kappa} &= \bm{K}, \\
\end{aligned} \qquad
\begin{aligned}
\bm\alpha_{\theta} &= \frac{\rho b^3}{12} \diffp{\bm{\theta}}{t}, \\
\bm\alpha_{\gamma} &= \bm{\gamma}. \\
\end{aligned}
\end{equation}
The co-energy variables are found by computing the variational derivatives of the Hamiltonian
\begin{equation}
\begin{aligned}
e_w &:= \diffd{H}{\alpha_w} = \diffp{w}{t},  \\
\bm{E}_{\kappa} &:= \diffd{H}{\bm{A}_{\kappa}} = \bm{M}, \\
\end{aligned} \qquad
\begin{aligned}
\bm{e}_{\theta} &:= \diffd{H}{\bm\alpha_{\theta}} = \diffp{\bm{\theta}}{t}, \\
\bm{e}_{\gamma} &:= \diffd{H}{\bm{\alpha}_{\bm{\gamma}}} = \bm{q}. \\
\end{aligned}
\end{equation}
Energy and co-energy variables are related by a positive symmetric operator $\bm{\alpha} = \mathcal{Q} \bm{e}$
\begin{equation}
\mathcal{Q} = \mathrm{diag}[(\rho b)^{-1}, \; (\rho b^3/12)^{-1} , \; \mathcal{D}, \; \mathcal{C}].
\end{equation}

The port-Hamiltonian system is expressed as follows 
\begin{equation}
\label{eq:PH_sys_Min_Ten}
\diffp{}{t}
\begin{pmatrix}
\alpha_w \\
\bm\alpha_\theta \\
\bm{A}_\kappa \\
\bm\alpha_{\gamma} \\
\end{pmatrix} = 
\underbrace{\begin{bmatrix}
	0  & 0  & 0  & \mathrm{div} \\
	0 & 0 &  \mathrm{Div} & \bm{I}_{2 \times 2}\\
	0  & \mathrm{Grad}  & 0  & 0\\
	\mathrm{grad} & -\bm{I}_{2 \times 2} &  0 & 0  \\
	\end{bmatrix}}_{\mathcal{J}}
\begin{pmatrix}
e_w \\
\bm{e}_{\theta} \\
\bm{E}_{\kappa} \\
\bm{e}_{\gamma} \\
\end{pmatrix} + 
\begin{pmatrix}
f \\
\bm{\tau} \\
0 \\
0 \\
\end{pmatrix}.
\end{equation}
\begin{remark}\label{rmk:stdir}
The force and torque $f, \bm{\tau}$ define a distributed control $\bm{u}_d$. Together with the collocated distributed output $\bm{y}_d=[e_w, \ \bm{e}_\theta]^\top$, this system defines a Dirac structure. By computing the power balance, it would be possible to add boundary variables $\bm{u}_\partial,  \bm{y}_\partial$, and, consequently, to obtain a Stokes-Dirac structure. However, in this paper we focus on clamped boundary condition, i.e.
\[
e_w|_{\partial \Omega} = 0, \quad \bm{e}_{\theta}|_{\partial \Omega} = 0 \implies \bm{u}_\partial = [e_w|_{\partial \Omega}, \ \bm{e}_\theta|_{\partial \Omega}]^\top \equiv 0.
\]
More general boundary conditions may be treated as well.
\end{remark}

\subsection{Kirchhoff plate}
The Kirchhoff plate model is a generalization to the 2D case of the Euler-Bernoulli beam model. The classical equations for this model are (\cite{timoshenko1959theory}) 
\begin{equation}
\label{eq:clKir}
\displaystyle \rho b \diffp[2]{w}{t} = -\mathrm{div}(\mathrm{Div}(\bm{M})) + f, \quad (\bm{x}, t) \in \Omega \times [0, t_f].
\end{equation}
As in the Mindlin model, the bending moment tensor and the curvature are related $\bm{M} = \mathcal{D} \bm{K} \in \mathbb{S}$ (with $\mathcal{D}$ defined in \eqref{eq:bend_rig_tensor}). Following the Kirchhoff assumption, the curvature tensor is the Hessian of the vertical displacement
\begin{equation*}
\bm{K} := \mathrm{Grad}(\mathrm{grad}(w)) \in \mathbb{S}.
\end{equation*}
 The kinetic and potential energy $E_c, E_p$ read
\begin{equation}
E_c =  \frac{1}{2} \int_{\Omega} \rho b \left(\diffp{w}{t} \right)^2 \d\Omega, \quad
E_p = \frac{1}{2} \int_{\Omega} \bm{M} \cddot \bm{K} \d\Omega.
\end{equation} 
The Hamiltonian is then given by $H=E_c + E_p$. Selecting as energy variables
\begin{equation}
\alpha_w = \rho b \diffp{w}{t}, \quad \bm{A}_{\kappa} = \bm{K}, 
\end{equation}
the co-energy variables are found by computing the variational derivatives of the Hamiltonian
\begin{equation}
e_w := \diffd{H}{\alpha_w} = \diffp{w}{t}, \quad \bm{E}_{\kappa} := \diffd{H}{\bm{A}_{\kappa}} = \bm{M}. \\
\end{equation}
The coercive operator linking energy and co-energy variables reads
\begin{equation}
\mathcal{Q} = \mathrm{diag}[(\rho b)^{-1}, \mathcal{D}].
\end{equation}
 
The port-Hamiltonian system is expressed as follows 
\begin{equation}
\label{eq:PH_sys_Kir_Ten}
\diffp{}{t}
\begin{pmatrix}
\alpha_w \\
\bm{A}_\kappa \\
\end{pmatrix} = 
\underbrace{\begin{bmatrix}
	0  & -\mathrm{div} \circ \mathrm{Div} \\
	\mathrm{Grad} \circ \mathrm{grad}  & 0 \\
	\end{bmatrix}}_{\mathcal{J}}
\begin{pmatrix}
e_w \\
\bm{E}_{\kappa} \\
\end{pmatrix}+ 
\begin{pmatrix}
f \\
0 \\
\end{pmatrix}.
\end{equation}
Following Remark \ref{rmk:stdir} this system would define a Stokes-Dirac if appropriate boundary variables were added. However, in this paper simply supported boundary conditions are considered, i.e.
\[
e_w|_{\partial \Omega} = 0, \quad m_{\text{nn}}|_{\partial \Omega}:= \bm{n}^\top \bm{E}_\kappa \bm{n}|_{\partial \Omega} = 0,
\]
hence no boundary control is present. Differently from the Mindlin plate case, generic boundary conditions demand an accurate analysis, see for instance \cite{Blum1990,mixed_kirchhoff}.

\section{Available mixed finite elements}
\label{sec:mixed}

In this section suitable semi-discretized models are derived. For the Mindlin model, two different formulation are presented: the first one enforces the symmetry of the momenta tensor strongly (\S\ref{sec:min_strong}), the second weakly (\S\ref{sec:min_weak}). For the Kirchhoff plate, the formulation is based on the the non-conforming Hellan-Herrmann-Johnson method (HHJ) (\S\ref{sec:HHJ}). 

\begin{remark}
System \eqref{eq:PH_sys_Min_Ten}, \eqref{eq:PH_sys_Kir_Ten} can be expressed using either the energy or the co-energy variables. The  most adapted formulation to the existing mixed finite element literature is the co-energy one, which reads $\mathcal{Q}^{-1} \partial_t \bm{e} = \mathcal{J} \bm{e}$.
\end{remark}

\subsection{Mindlin plate with strongly imposed symmetry}\label{sec:min_strong}

The weak formulation with strongly imposed symmetry seeks $\{e_w, \bm{e}_{\bm{\theta}}, \bm{E}_{\kappa}, \bm{e}_{\gamma}\} \in L \times V \times S \times D$ so that 
\begin{equation}
\label{eq:weak_min_PH_strong}
\begin{aligned}
(v_w, \ \rho b \dot{e}_w) &= (v_w, \mathrm{div} \bm{e}_\gamma) + (v_w, f), \\ 
(\bm{v}_\theta, \ \rho b^3/12  \dot{\bm{e}}_\theta) &= (\bm{v}_\theta, \mathrm{Div} \bm{E}_\kappa + \bm{e}_\gamma) + (\bm{v}_\theta, \bm{\tau}), \\  
(\bm{V}_\kappa, \ \mathcal{D}^{-1} \dot{\bm{E}}_\kappa) &= -(\mathrm{Div} \bm{V}_\kappa,  \bm{e}_\theta), \\ 
(\bm{v}_\gamma, \ \mathcal{C}^{-1} \dot{\bm{e}}_\gamma) &= -(\mathrm{div} \bm{v}_\gamma, e_w ) + (\bm{v}_\gamma, \bm{e}_{\theta}), \\ 
\end{aligned} \quad
\begin{aligned}
v_w \in L, \\
\bm{v}_\theta \in V, \\
\bm{V}_\kappa \in S, \\
\bm{v}_\gamma \in D.
\end{aligned}
\end{equation}
This formulation is obtained by multiplying each equation by a test function belonging to the same space as the corresponding unknown, and integrating over the domain. The final system is obtained by integrating by parts the last two lines of \eqref{eq:PH_sys_Min_Ten} and considering clamped boundary conditions. Obtaining stable finite elements that embed the symmetry of the stress tensor for the elastodynamics problem has proven to be a difficult task. The easiest implementation is the one presented in \cite{becacheWave,becacheElas}. The main disadvantage is that this scheme requires the domain to be given by a union of rectangles, as the mesh elements have to be square. However, this allows constructing a simple element for the momenta tensor. The polynomial spaces for the discretization are
\[
N_{k} = \{p(x, y)| \; p(x, y) = \sum_{i\le k, j\le k} a_{ij} x^i y^j  \}.
\]
Given a regular mesh $\mathcal{R}_h$ with square elements $Q$ the following spaces are introduced as discretization spaces
\begin{equation}
\label{eq:BTJ}
	\begin{aligned}
	L_h^{\text{BJT}} &= \{w_h \in L | \ \forall Q, \ w_h|_{Q} \in N_{k-1} \}, \\
	V_h^{\text{BJT}} &= \{\bm{\theta}_h \in V | \ \forall Q,\ \bm{\theta}_h|_{Q} \in (N_{k-1})^2 \}, \\
	S_h^{\text{BJT}} &= \{m_{12} \in H^1(\Omega)| \ \forall Q,\ m_{12}|_{Q} \in N_{k} \}  \\
	&\cup \{(m_{11}, m_{22}) \in D| \; \forall Q,\ (m_{11}, m_{22})|_{Q} \in N_{k} \}, \\
	D_h^{\text{BJT}} &= \{\bm{q}_h \in D | \ \forall Q,\ \bm{q}_h|_{Q} \in N_{k} \}, \\ 
	\end{aligned}
\end{equation}
where BTJ stands for the initials of the authors in \cite{becacheWave,becacheElas}. Combining the results of both papers, the following error estimates are conjectured:
\begin{conjecture}
Assuming a smooth solution to problem~\eqref{eq:weak_min_PH_strong}, the following error estimates hold 
\begin{equation}
\label{eq:errBEC}
\begin{aligned}
||e_w - e_w^h||_{L^{\infty}(L^2)} &\lesssim h^{k}, \\
||\bm{e}_\theta - \bm{e}_\theta^h||_{L^{\infty}(L^2)} &\lesssim h^{k}, \\
\end{aligned} \quad
\begin{aligned}
||\bm{E}_\kappa - \bm{E}_\kappa^h||_{L^{\infty}(L^2)} &\lesssim  h^{k}, \\
||\bm{e}_\gamma - \bm{e}_\gamma^ h||_{L^{\infty}(L^2)} &\lesssim  h^{k}, \\
\end{aligned} 
\end{equation}
where the notation $a \lesssim  b$ means $a \le C b$. The constant $C$ depends only on the true solution and on the final time.
\end{conjecture}

\subsection{Mindlin plate with weakly imposed symmetry}\label{sec:min_weak}
Formulation \eqref{eq:weak_min_PH_strong} has to be modified to impose the symmetry of the momenta tensor weakly. Taking the weak form of the third equation in \eqref{eq:PH_sys_Min_Ten}, we get
\[
(\bm{V}_\kappa, \ \mathcal{D}^{-1} \dot{\bm{E}}_\kappa) = (\bm{V}_\kappa, \mathrm{Grad}\ \bm{e}_\theta ). 
\]
The symmetric gradient can be rewritten as 
\[
\mathrm{Grad}\ \bm{\theta} = \mathrm{grad}\ \bm{\theta} - \mathrm{skw}(\mathrm{grad} \ \bm{\theta}),
\]
where $\mathrm{skw}(\bm{A})=(\bm{A} - \bm{A}^\top)/2$ is the skew-symmetric part of matrix $\bm{A}$. Introducing the new variable $\bm{E}_r = \mathrm{skw}(\mathrm{grad}\ \bm{\theta})$, then $\{\bm{e}_\theta, \bm{E}_\kappa, \bm{E}_r\} \in V\times M\times K$ satisfy (remind that $\bm{e}_\theta = \partial_t {\bm{\theta}}$)
\begin{equation*}
	\begin{aligned}
	(\bm{V}_\kappa, \ \mathcal{D}^{-1} \dot{\bm{E}}_\kappa) &= (\bm{V}_\kappa, \mathrm{grad}\ \bm{e}_\theta) - (\bm{V}_\kappa, \ \dot{\bm{E}}_r), \\
	&= -(\mathrm{Div}\bm{V}_\kappa, \bm{e}_\theta) - (\bm{V}_\kappa, \ \dot{\bm{E}}_r).
	\end{aligned}
\end{equation*}
The momenta tensor is weakly symmetric if $(\bm{V}_r, \ \bm{E}_{\kappa}) = 0$. The weak formulation then consists in finding \\
$\{e_w, \bm{e}_{\bm{\theta}}, \bm{E}_{\kappa}, \bm{e}_{\gamma}, \bm{E}_{r}\}$ in $L \times V \times M \times D \times K$ so that 
 \begin{equation}
 \label{eq:weak_min_PH_weak}
 \begin{aligned}
 (v_w, \ \rho b \dot{e}_w) &= (v_w, \mathrm{div} \bm{e}_\gamma) + (v_w, f), \\ 
 (\bm{v}_\theta, \ \rho b^3/12  \dot{\bm{e}}_\theta) &= (\bm{v}_\theta, \mathrm{Div} \bm{E}_\kappa + \bm{e}_\gamma) + (\bm{v}_\theta, \bm{\tau}), \\  
 (\bm{V}_\kappa, \ \mathcal{D}^{-1} \dot{\bm{E}}_\kappa) &= -(\mathrm{Div} \bm{V}_\kappa,  \bm{e}_\theta) - (\bm{V}_\kappa, \ \dot{\bm{E}}_r), \\ 
 (\bm{v}_\gamma, \ \mathcal{C}^{-1} \dot{\bm{e}}_\gamma) &= -(\mathrm{div} \bm{v}_\gamma, e_w ) + (\bm{v}_\gamma, \bm{e}_{\theta}), \\ 
 (\bm{V}_r, \ \dot{\bm{E}}_\kappa) &= 0
 \end{aligned} \quad
 \begin{aligned}
 v_w \in L, \\
 \bm{v}_\theta \in V, \\
 \bm{V}_\kappa \in S, \\
 \bm{v}_\gamma \in D, \\
 \bm{V}_r \in K, \\
 \end{aligned}
 \end{equation}
Consider a regular triangulation $\mathcal{T}_h$ with elements $T$. The space of polynomials of order $k$ on a mesh cell is denoted by $P_k$. The following spaces are used as discretization spaces
\begin{equation}
\label{eq:AFW}
\begin{aligned}
L_h^{\text{AFW}} &= \{w_h \in L | \ \forall T, \ w_h|_{T} \in P_{k-1} \}, \\
V_h^{\text{AFW}} &= \{\bm{\theta}_h \in V | \ \forall T,\ \bm{\theta}_h|_{T} \in (P_{k-1})^2 \}, \\
S_h^{\text{AFW}} &= \{(m_{11}, m_{12}) \in D| \ \forall T,\ (m_{11}, m_{12})|_{T} \in BDM_{[k]} \}  \\
& \cup \{(m_{21}, m_{22}) \in D| \ \forall T,\ (m_{21}, m_{22})|_{T} \in BDM_{[k]} \}, \\
D_h^{\text{AFW}} &= \{\bm{q}_h \in D | \ \forall T,\ \bm{q}_h|_{T} \in RT_{[k-1]} \}, \\
K_h^{\text{AFW}} &= \{\bm{R}_h \in K | \ \forall T, \ w_h|_{T} \in P_{k-1} \}, \\ 
\end{aligned}
\end{equation}
where $BDM$ are the Brezzi-Douglas-Marini elements and $RT$ the Raviart-Thomas elements. The acronym AFW stands for Arnold-Falk-Winther. A convergence analysis for the general elastodynamics problem with weak symmetry in the $L^\infty (L^2)$ norm is detailed \cite{ArnoldWeak}. A convergence study for the wave equation with mixed finite elements in the $L^\infty (L^2)$ is presented in \cite{Geveci}. Combining the results of the two, the following error estimates are conjectured:
\begin{conjecture}
	Assuming a smooth solution to problem~\eqref{eq:weak_min_PH_strong}, the following error estimates hold 
	\begin{equation}
	\label{eq:errAFW}
	\begin{aligned}
	||e_w - e_w^h||_{L^{\infty}(L^2)} &\lesssim h^{k}, \\
	||\bm{e}_\theta - \bm{e}_\theta^h||_{L^{\infty}(L^2)} &\lesssim h^{k}, \\
	||\bm{E}_r - \bm{E}_r^h||_{L^{\infty}(L^2)} &\lesssim h^{k}. \\
	\end{aligned} \quad
	\begin{aligned}
	||\bm{E}_\kappa - \bm{E}_\kappa^h||_{L^{\infty}(L^2)} &\lesssim  h^{k}, \\
	||\bm{e}_\gamma - \bm{e}_\gamma^ h||_{L^{\infty}(L^2)} &\lesssim  h^{k}, \\
	\end{aligned} 
	\end{equation}
\end{conjecture}

\subsection{The HHJ scheme for the Kirchhoff plate}\label{sec:HHJ}
For the Kirchhoff plate, the HHJ scheme can be used to obtain a structure-preserving discretization. The discussion follows \cite{arnold2019hellan}. Given the non conforming nature of this scheme, it is necessary to first introduce the discrete functional spaces and state the problem directly in discrete form. The vertical displacement is approximated using continuous Lagrange polynomials, while the momenta tensor is discretized using the HHJ element
\begin{equation}
\label{eq:HHJ}
\begin{aligned}
W_h = &\{w_h \in H^1_0(\Omega)| \; \forall T, \; w_h|_{T} \in P_{k} \}, \\
U_h = &\{\bm{M}_h \in L^2(\Omega, \mathbb{S})| \; \forall T, \; \bm{M}_h|_{T} \in P_{k-1}(\mathbb{S}) , \\ 
&\, \ \bm{M}_h \text{ is normal-normal continous across elements}\}.
\end{aligned}
\end{equation}
The normal to normal continuity means that if two triangles $T_1, T_2$ share a common edge $E$ then $\bm{n}^\top (\bm{M}_h|_{T_1}) \bm{n} = \bm{n}^\top (\bm{M}_h|_{T_2}) \bm{n}$ on $E$. Taking system \eqref{eq:PH_sys_Kir_Ten} and multiplying the first equation by $v_w \in W_h$ and integrating over a triangle
\begin{equation*}
	\begin{aligned}
	& - (v_w, \ \mathrm{div}\mathrm{Div} \bm{E}_\kappa))_{T} = (\nabla v_w, \ \mathrm{Div} \bm{E}_\kappa))_{T}=, \\
	& -(\nabla^2 v_w, \ \bm{E}_\kappa)_T + (\partial_n v_w, \bm{n}^\top\bm{E}_\kappa \bm{n})_{\partial T} + (\partial_s v_w, \bm{s}^\top\bm{E}_\kappa \bm{n})_{\partial T}. \\
	\end{aligned}
\end{equation*}
A double integration by parts is applied to get the final equation. Summing up over all triangles provides for the penultimate term
\begin{equation*}
\sum_{T \in \mathcal{T}_h} (\partial_n v_w, \bm{n}^\top\bm{E}_\kappa \bm{n})_{\partial T} = \sum_{E \in \mathcal{E}_h} ([\![\partial_n v_w]\!], m_{\text{nn}})_{E},
\end{equation*} 
where $\mathcal{E}_h$ is the set of all edges belonging to the mesh and $[\![a]\!] = a|_{T_1} + a|_{T_2}$ denotes the jump of a function across shared edges. For a boundary edge it is simply the value of the function. For the final term, it holds $(\partial_s v_w, \bm{s}^\top\bm{E}_\kappa \bm{n})_{\partial T}=0$, since $v_w$ is continuous across the edge boundaries and the normal switches sign. We are now in a position to state the final weak form. Given the definition
\[
b_h(v_w, \ \bm{E}_{\kappa}) := - \sum_{T \in \mathcal{T}_h} ( \nabla^2 v_w, \ \bm{E}_\kappa) + \sum_{E \in \mathcal{E}_h} ([\![\partial_n v_w]\!], m_{\text{nn}})_{E}, 
\]
find $(e_w, \bm{E}_\kappa) \in W_h \times U_h$ such that
\begin{equation}
\label{eq:weak_kir_PH}
\begin{aligned}
(v_w, \ \rho b \dot{e}_w) &= +b_h(v_w, \ \bm{E}_{\kappa}) + (v_w, f), \\ 
(\bm{V}_\kappa, \ \mathcal{D}^{-1} \dot{\bm{E}}_\kappa) &= -b_h(e_w, \ \bm{V}_{\kappa}), \\ 
\end{aligned} \quad
\begin{aligned}
v_w \in W_h, \\
\bm{V}_\kappa \in U_h. \\
\end{aligned}
\end{equation}

For the associated static problem, under the hypothesis of smooth solutions, optimal convergence of order $O(k)$ for $w \in H^1$ and $\bm{M} \in L^2$ has been established. So, it is natural to conjecture the following result for the dynamic problem:
\begin{conjecture}
Assuming a smooth solution for problem~\eqref{eq:weak_kir_PH}, the following error estimates hold
\begin{equation}
\label{eq:errHHJ}
||e_w - e_w^h||_{L^{\infty} (H^1)} \lesssim h^{k}, \qquad
||\bm{E}_\kappa - \bm{E}_\kappa^h||_{L^{\infty} (L^2)} \lesssim h^{k}.
\end{equation}
\end{conjecture}

\section{Numerical experiments}
\label{sec:numerics}
In this section numerical test cases are used to verify the conjectured orders of convergence for the two problems. Upon discretization, system \eqref{eq:weak_min_PH_strong}, \eqref{eq:weak_min_PH_weak}, \eqref{eq:weak_kir_PH} assumes the form 
\[
M \dot{\bm{e}} = J \bm{e}.
\]
Matrix $J$ is skew-symmetric, matrix $M$ is symmetric and positive definite for \eqref{eq:weak_min_PH_strong}, \eqref{eq:weak_kir_PH} while it is symmetric but indefinite for \eqref{eq:weak_min_PH_weak}, because of the multiplier that enforces the symmetry. The Firedrake library (\cite{rathgeber2017firedrake}) is used to generate the matrices. To integrate the equations in time a Crank-Nicholson scheme has been used, for all simulations. The time step is set to $\Delta t = h/10$ to have a lower impact of the time discretization error with respect to the spatial error. The final time is set to one $t_f = 1 [\textrm{s}]$ for all simulations. To compute the $L^\infty (\mathcal{X})$ space-time dependent norm  the discrete norm $L^\infty_{\Delta t} (\mathcal{X})$ is used
\[
||\cdot ||_{L^\infty (\mathcal{X})} \approx || \cdot ||_{L^\infty_{\Delta t} (\mathcal{X})} = \max_{t \in t_i} ||\cdot||_{\mathcal{X}},
\]
where $t_i$ are the discrete simulation instants. 
\subsection{Numerical test for the Mindlin plate}
Constructing an analytical solution for a vibrating Mindlin plate is far from trivial. Therefore, the solution for the static case presented in \cite{mindlinVeiga} is exploited. \\
\textbf{Step 1 } Consider a distributed static force given by 
\begin{equation*}
\begin{aligned}
f_s(x,y)=\frac{E_Y}{12 (1-\nu^2)} \{12 y(y-1)(5x^2-5x+1)\\
\times [2y^2(y-1)2+x(x-1)(5y^2-5y+1)] +12x(x-1)\\
\times (5y^2-5y+1)[2x^2(x-1)2+y(y-1)(5x^2-5x+1)]\}.
\end{aligned}
\end{equation*}
The static displacement and rotation are given by
\begin{align*}
	w_s(x,y) &= \frac{1}{3} x^3(x-1)^3 y^3 (y-1)^3\\
	&-\frac{2 b^2}{5(1-\nu)}[y^3(y-1)^3 x(x-1)(5 x^2-5x+1). \\
	\bm{\theta}_{s}(x,y) &= 
	\begin{pmatrix}
	y^3(y-1)^3 \ x^2 (x-1)^2 (2x-1) \\
	x^3(x-1)^3 \ y^2 (y-1)^2 (2y-1) \\
	\end{pmatrix}
\end{align*}
The static solution solves the following problem defined on the square domain $\Omega=(0,1)\times(0,1)$:
\begin{equation}
\begin{aligned}
0 &= \mathrm{div} \ \bm{q}_s + f_s , \\
0 &= \mathrm{Div} \bm{M}_s + \bm{q}_s, \\
\end{aligned} \qquad
\begin{aligned}
\mathcal{D}^{-1} \bm{M}_s &= \mathrm{Grad} \ \bm{\theta}_s, \\
\mathcal{C}^{-1} \bm{q}_s &= \mathrm{grad} \ w_s - \bm{\theta}_s. \\
\end{aligned}
\end{equation}
\textbf{Step 2 } Given the linear nature of the system a solution for the dynamic problem is found by multiplying the static solution by a time dependent term. For simplicity a sinus function is chosen
\[
w_d(x,y,t) = w_s(x,y) \sin(t), \quad \bm{\theta}_d(x,y,t) = \bm\theta_s(x,y) \sin(t).
\]
For the port-Hamiltonian system velocities are needed
\[
e_w^\text{ex}(x,y,t) = w_s(x,y) \cos(t), \quad \bm{e}_\theta^\text{ex}(x,y,t) = \bm\theta_s(x,y) \cos(t).
\]
The momenta and shear force are then defined by
\[
\bm{M}_d = \bm{E}_\kappa^\text{ex} =  \mathcal{D} \ \mathrm{Grad} \ \bm{\theta}_d, \quad \bm{q}_d = \bm{e}_\gamma^\text{ex} = \mathcal{C}(\mathrm{grad} \ w_d - \bm{\theta}_d)
\]
\textbf{Step 3 } Appropriate forcing terms have to be introduced (i.e. $f, \bm{\tau}$ in \eqref{eq:clMin}). The force and torque in the dynamical case become
\begin{equation*}
f_d = f_s \sin(t) + \rho b \partial_{tt} w_d, \qquad
\bm{\tau}_d = \frac{\rho b^3}{12} \partial_{tt} \bm{\theta}_d.
\end{equation*}
Variables $(e_w^\text{ex}, \bm{e}_\theta^\text{ex}, \bm{E}_\kappa^\text{ex}, \bm{e}_\gamma^\text{ex})$ under excitations $(f_d, \bm{\tau}_d)$ solve problem~\eqref{eq:PH_sys_Min_Ten}. The solution being smooth, the conjectured error estimates should hold. The numerical values of the physical parameters are reported in Table \ref{tab:parMin}.

\begin{table}[h]
	\centering
	\begin{tabular}{ccccc}
		\hline 
		\multicolumn{5}{c}{Plate parameters} \\ 
		\hline 
		$E$ & $\rho$ & $\nu$ & $k$ & $h$ \\
		1 $[\textrm{Pa}]$ & $1\; [\textrm{kg}/\textrm{m}^3]$ & 0.3 & 5/6 & 0.1 $[\textrm{m}]$\\ 
		\hline 
	\end{tabular} 
	\captionsetup{width=0.95\linewidth}
	\vspace{1mm}
	\captionof{table}{Physical parameters for the Mindlin plate.}
	\label{tab:parMin}
\end{table}

\subsubsection{Results for the strong symmetry formulation} 

\begin{figure}[ht]%
	\centering
	\subfloat[][$L^\infty_{\Delta t} (L^2)$ error for $e_w$]{%
		\label{fig:errBEC1}%
		\includegraphics[width=0.48\columnwidth]{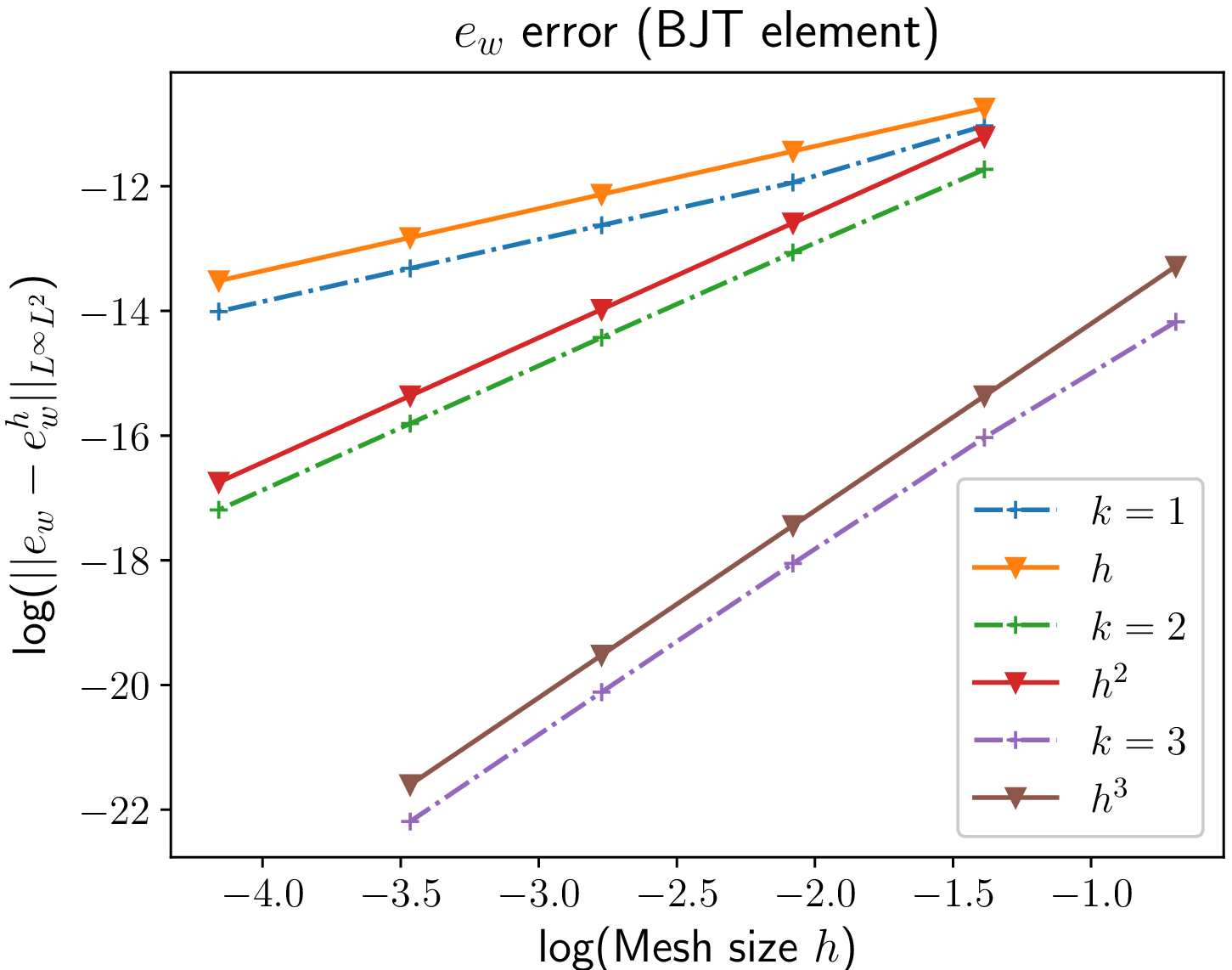}}%
	\hspace{8pt}%
	\subfloat[][$L^\infty_{\Delta t} (L^2)$ error for $\bm{e}_\theta$]{%
		\label{fig:errBEC2}%
		\includegraphics[width=0.48\columnwidth]{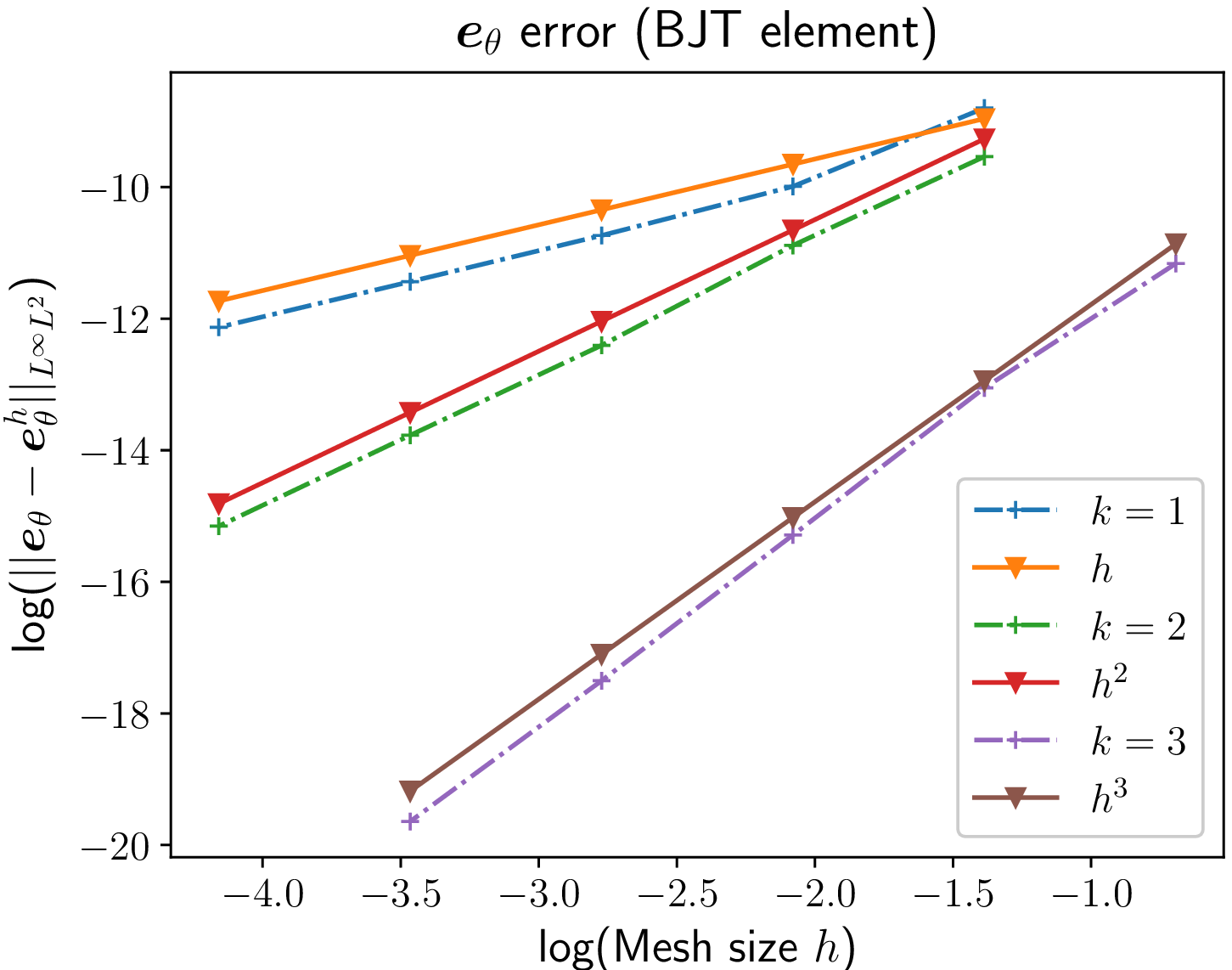}} \\
	\subfloat[][$L^\infty_{\Delta t} (L^2)$ error for $\bm{E}_\kappa$]{%
		\label{fig:errBEC3}%
		\includegraphics[width=0.48\columnwidth]{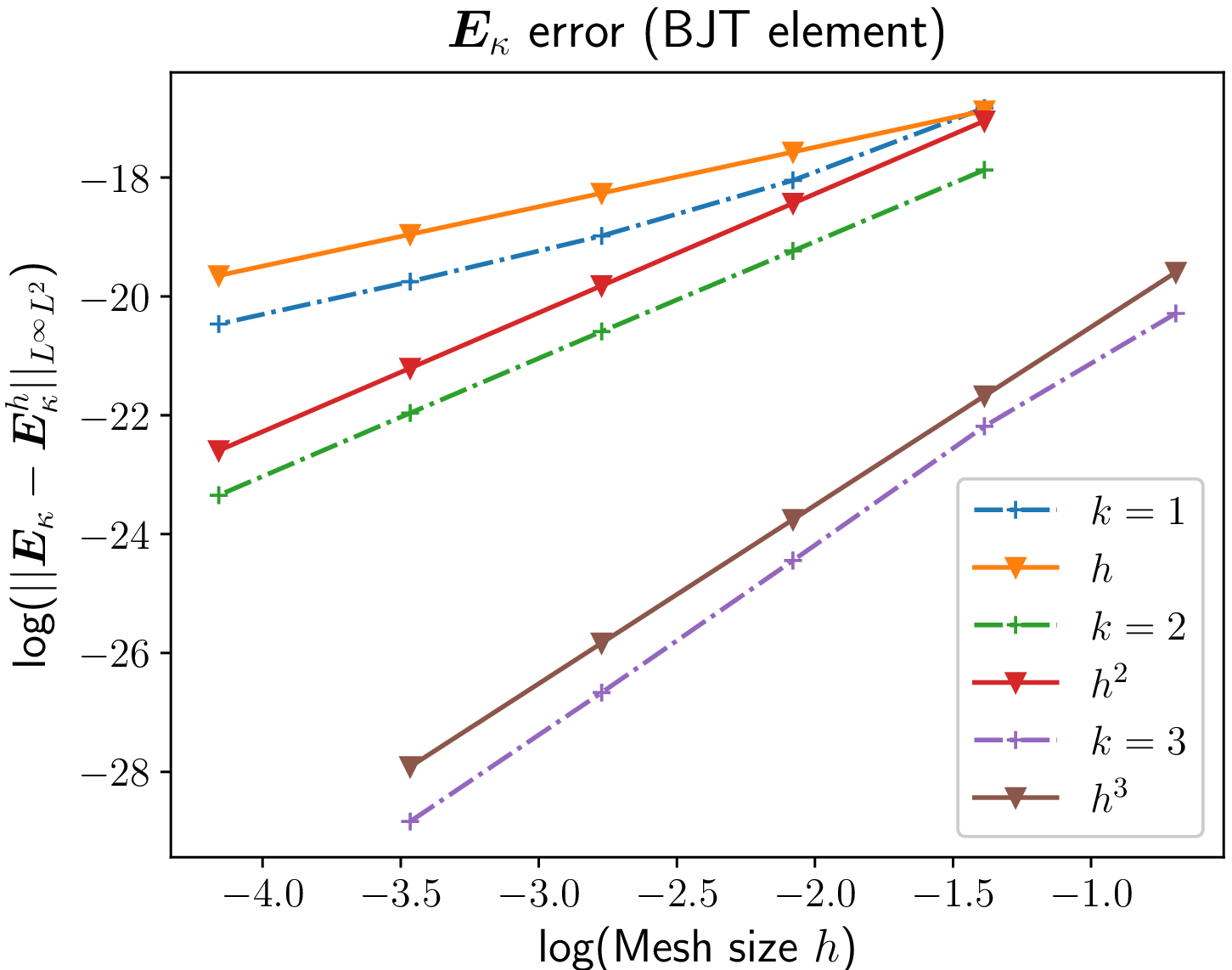}}%
	\hspace{8pt}%
	\subfloat[][$L^\infty_{\Delta t} (L^2)$ error for $\bm{e}_\gamma$]{%
		\label{fig:errBEC4}%
		\includegraphics[width=0.48\columnwidth]{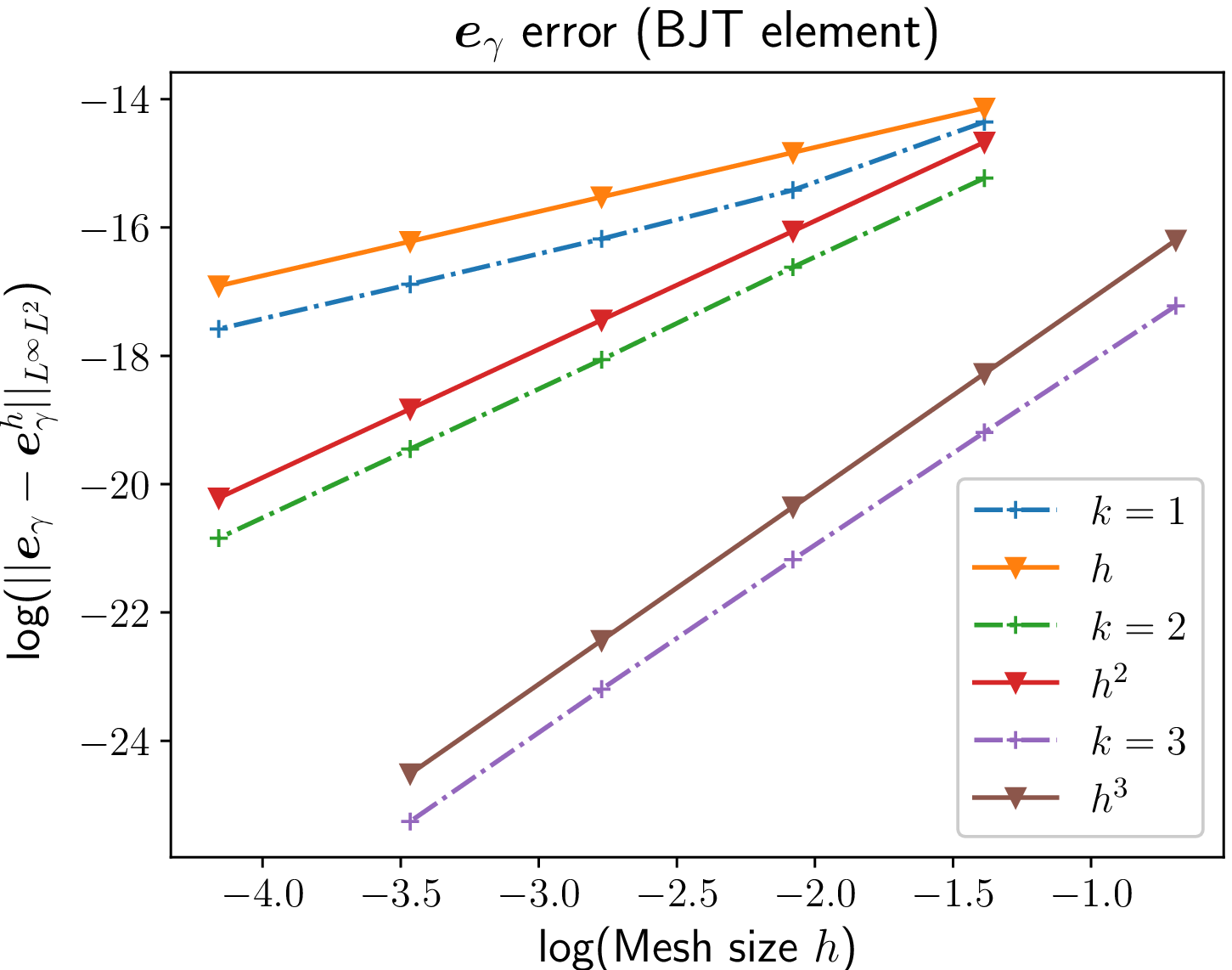}}%
	\caption[errorBEC]{Error for the Mindlin plate using the BJT elements}%
	\label{fig:errorBEC}%
\end{figure}

The weak form \eqref{eq:weak_min_PH_strong} and its corresponding finite elements \eqref{eq:BTJ} was implemented using Firedrake extruded mesh functionality (\cite{firedrake_extruded}). A direct solver based on an LU preconditioner is used. In Fig. \ref{fig:errorBEC} the errors for $(e_w, \bm{e}_\theta, \bm{E}_\kappa, \bm{e}_\gamma)$ are reported. As one can notice, the conjectured error estimates \eqref{eq:errBEC} are respected for all variables. 

\subsubsection{Results for the weak symmetry formulation} 
\begin{figure}[ht]%
	\centering
	\subfloat[][$L^\infty_{\Delta t} (L^2)$ error for $e_w$]{%
		\label{fig:errAFW1}%
		\includegraphics[width=0.48\columnwidth]{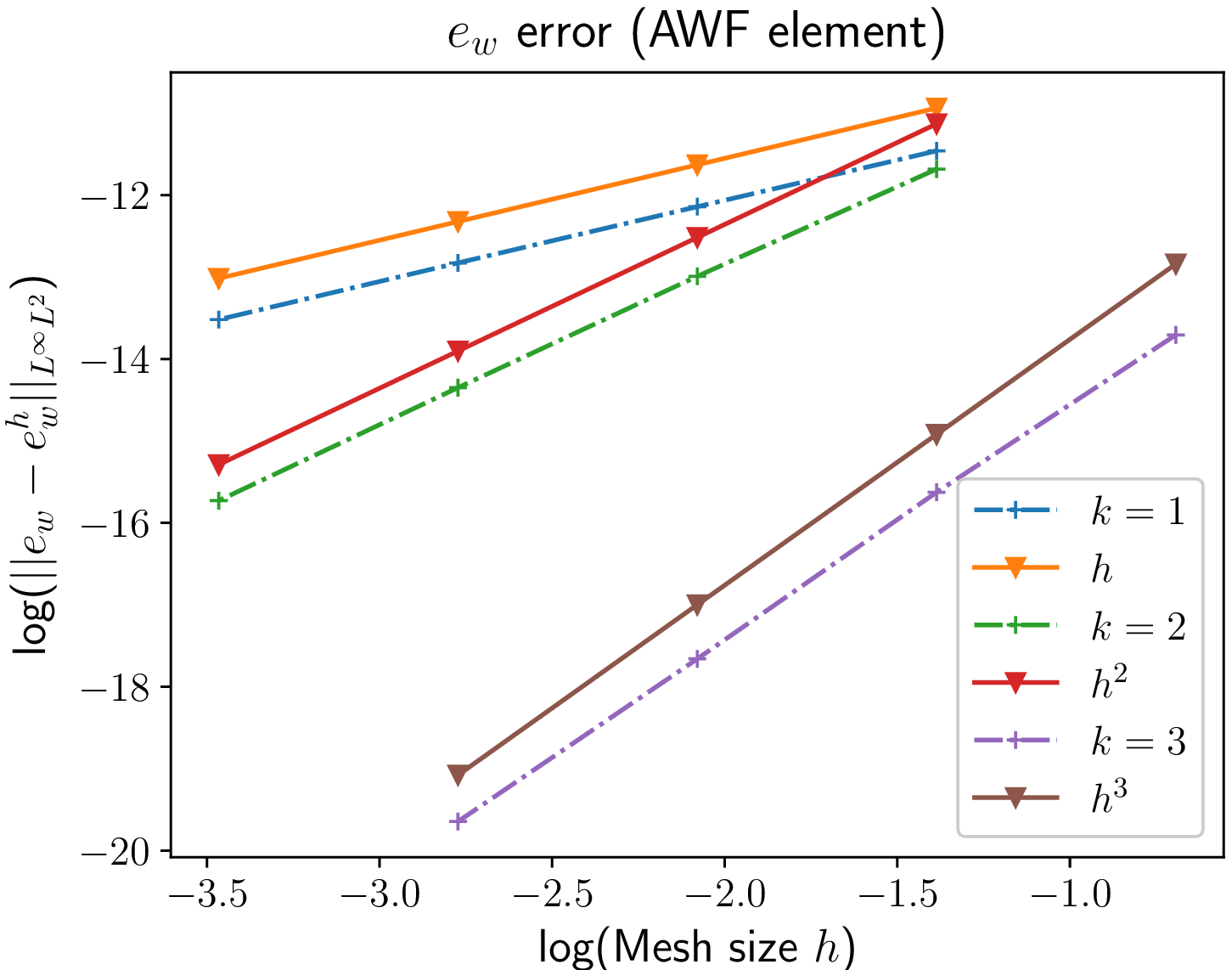}}%
	\hspace{8pt}%
	\subfloat[][$L^\infty_{\Delta t} (L^2)$ error for $\bm{e}_\theta$]{%
		\label{fig:errAFW2}%
		\includegraphics[width=0.48\columnwidth]{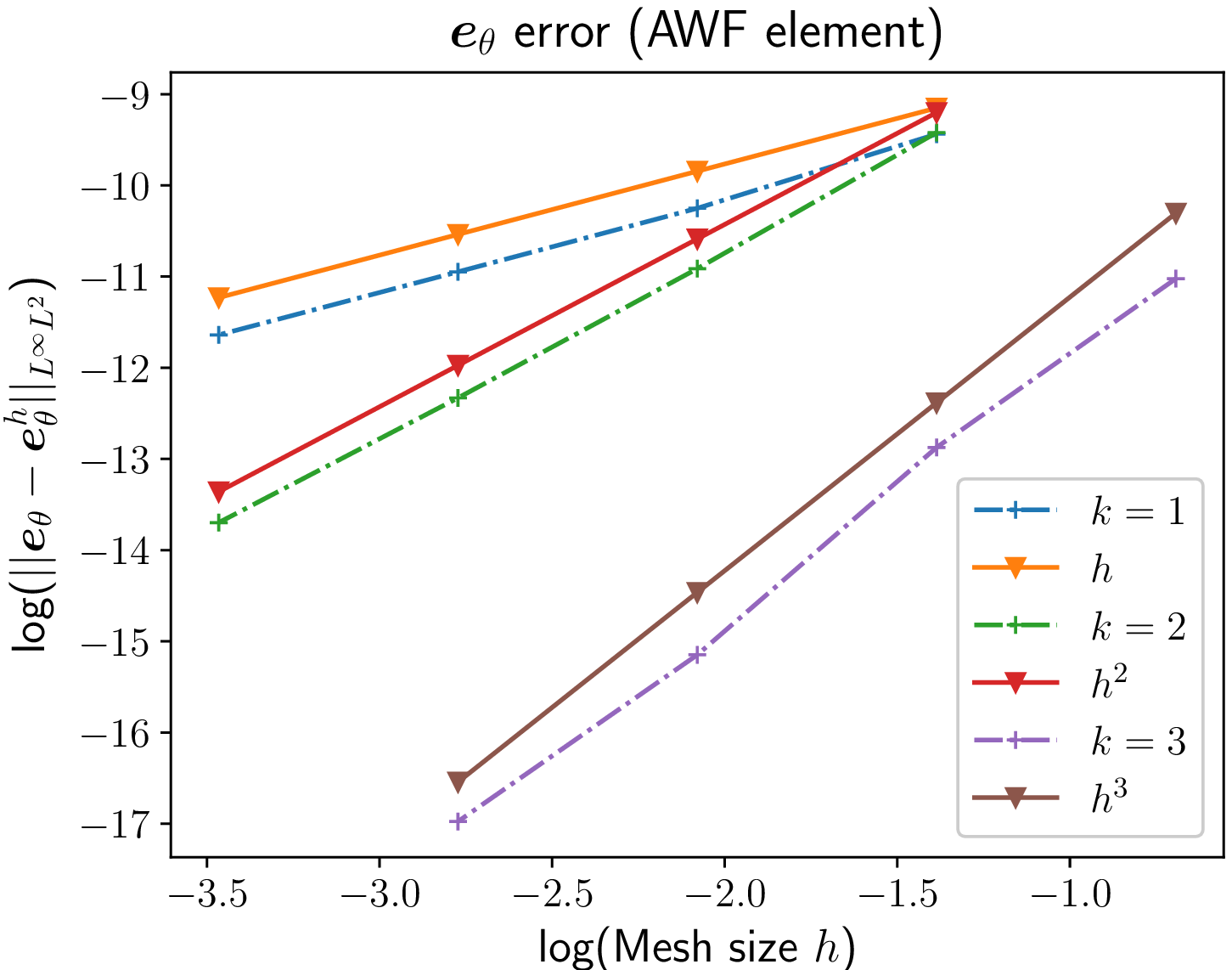}} \\
	\subfloat[][$L^\infty_{\Delta t} (L^2)$ error for $\bm{E}_\kappa$]{%
		\label{fig:errAFW3}%
		\includegraphics[width=0.48\columnwidth]{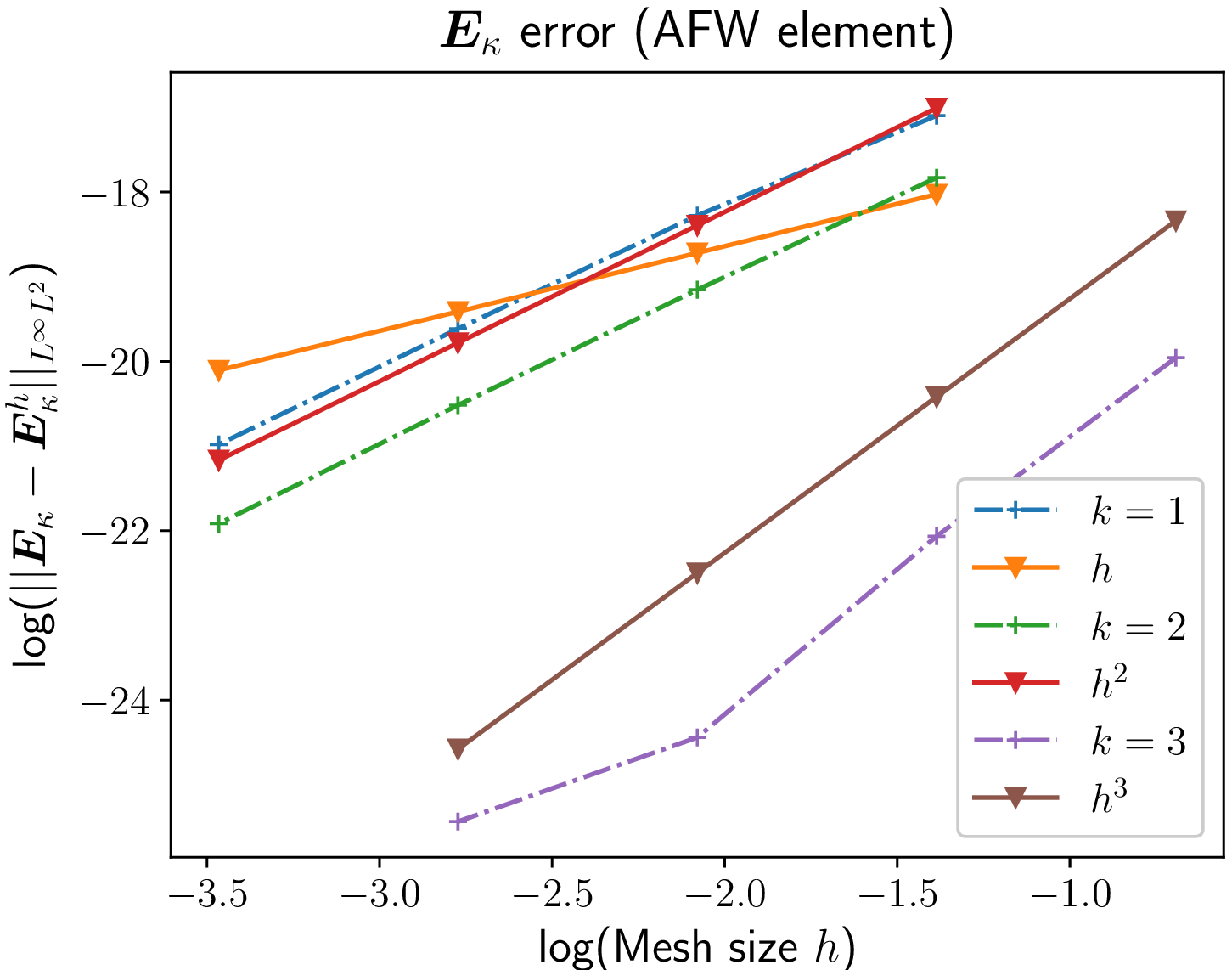}}%
	\hspace{8pt}%
	\subfloat[][$L^\infty_{\Delta t} (L^2)$ error for $\bm{e}_\gamma$]{%
		\label{fig:errAFW4}%
		\includegraphics[width=0.48\columnwidth]{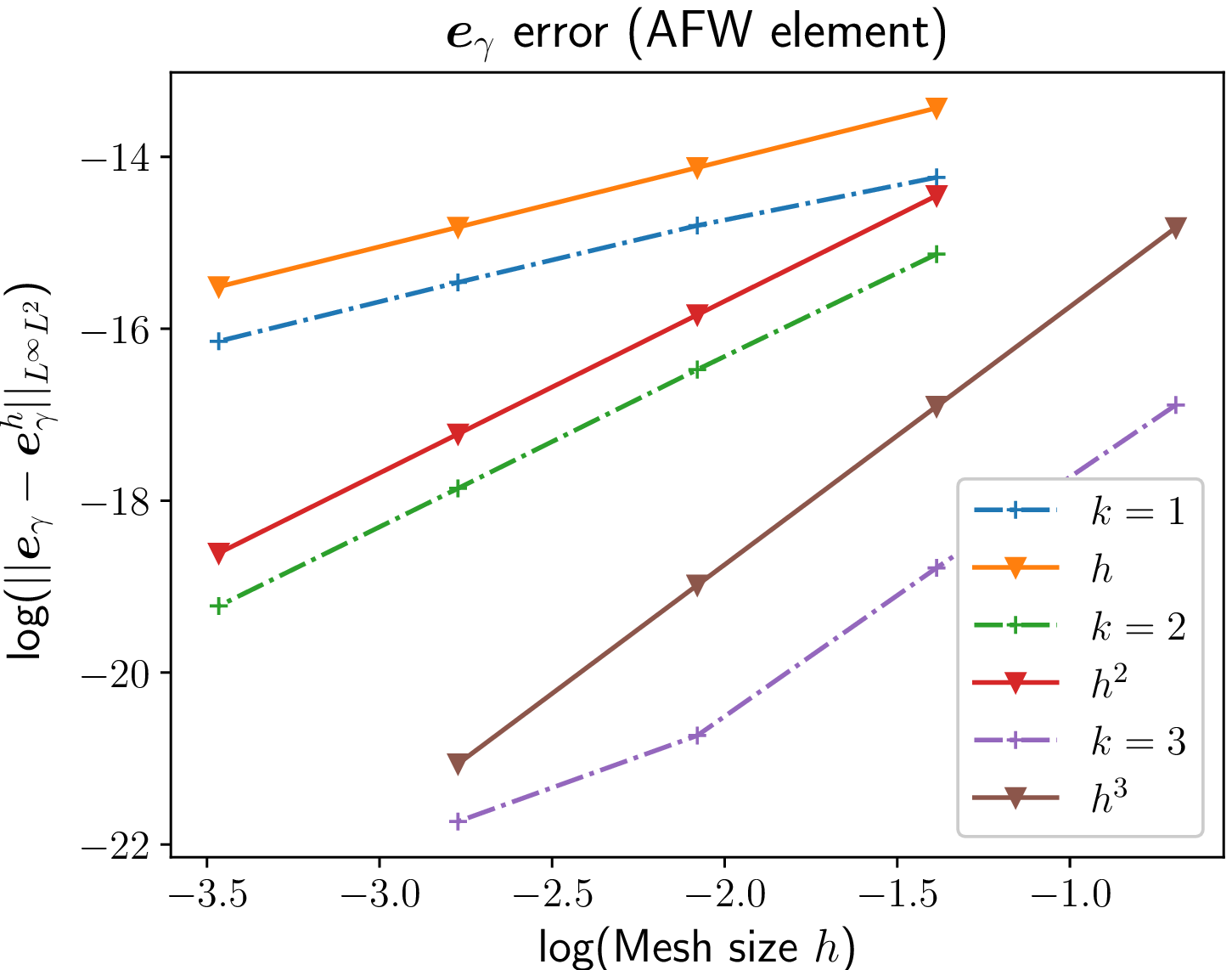}}%
	\caption[errorAFW]{Error for the Mindlin plate using the AFW elements}%
	\label{fig:errorAFW}%
\end{figure}

Formulation \eqref{eq:weak_min_PH_weak} and its element \eqref{eq:AFW} are considered here. A direct solver failed for high order cases (i.e. $k=3$). For this reason a generalized minimal residual method is used with restart number of iterations equal to 100. In Fig. \ref{fig:errorAFW} the errors for variables $(e_w, \bm{e}_\theta, \bm{E}_\kappa, \bm{e}_\gamma)$  are reported. The errors for $(e_w, \bm{e}_\theta, \bm{e}_\gamma)$ respect the conjectured result \eqref{eq:errAFW}. Variable $\bm{E}_\kappa$ exhibit a superconvergence phenomenon for the case $k=1$. In \cite{ArnoldWeak} no numerical study was carried out for the case $k=1$. The $BDM$ elements might be responsible for such superconvergence. The convergence order of $(\bm{E}_\kappa, \bm{e}_\gamma)$ deteriorates for $k=3$ for the finest mesh. This must be linked to errors due to the underlying large saddle-point problem. Indeed in \cite{ArnoldWeak} an hybridization method is used to transform the saddle-point problem into a symmetric positive definite one.

\subsection{Numerical test for the Kirchhoff plate}
\begin{figure}[ht]%
	\centering
	\subfloat[][$L^\infty_{\Delta t} (H^1)$ error for $e_w$]{%
		\label{fig:errHHJ1}%
		\includegraphics[width=0.48\columnwidth]{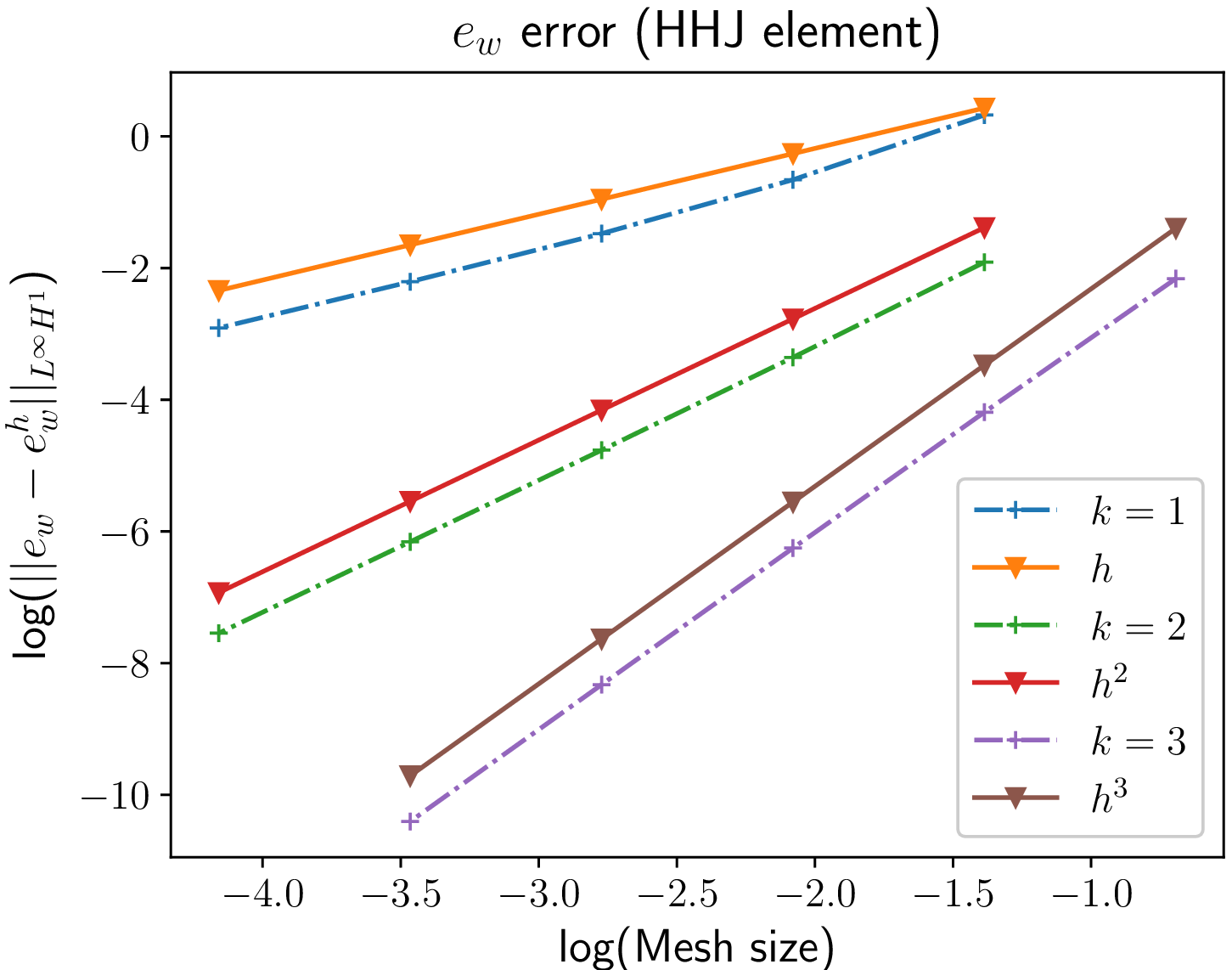}}%
	\hspace{8pt}%
	\subfloat[][$L^\infty_{\Delta t} (L^2)$ error for $\bm{E}_\kappa$]{%
		\label{fig:errHHJ2}%
		\includegraphics[width=0.48\columnwidth]{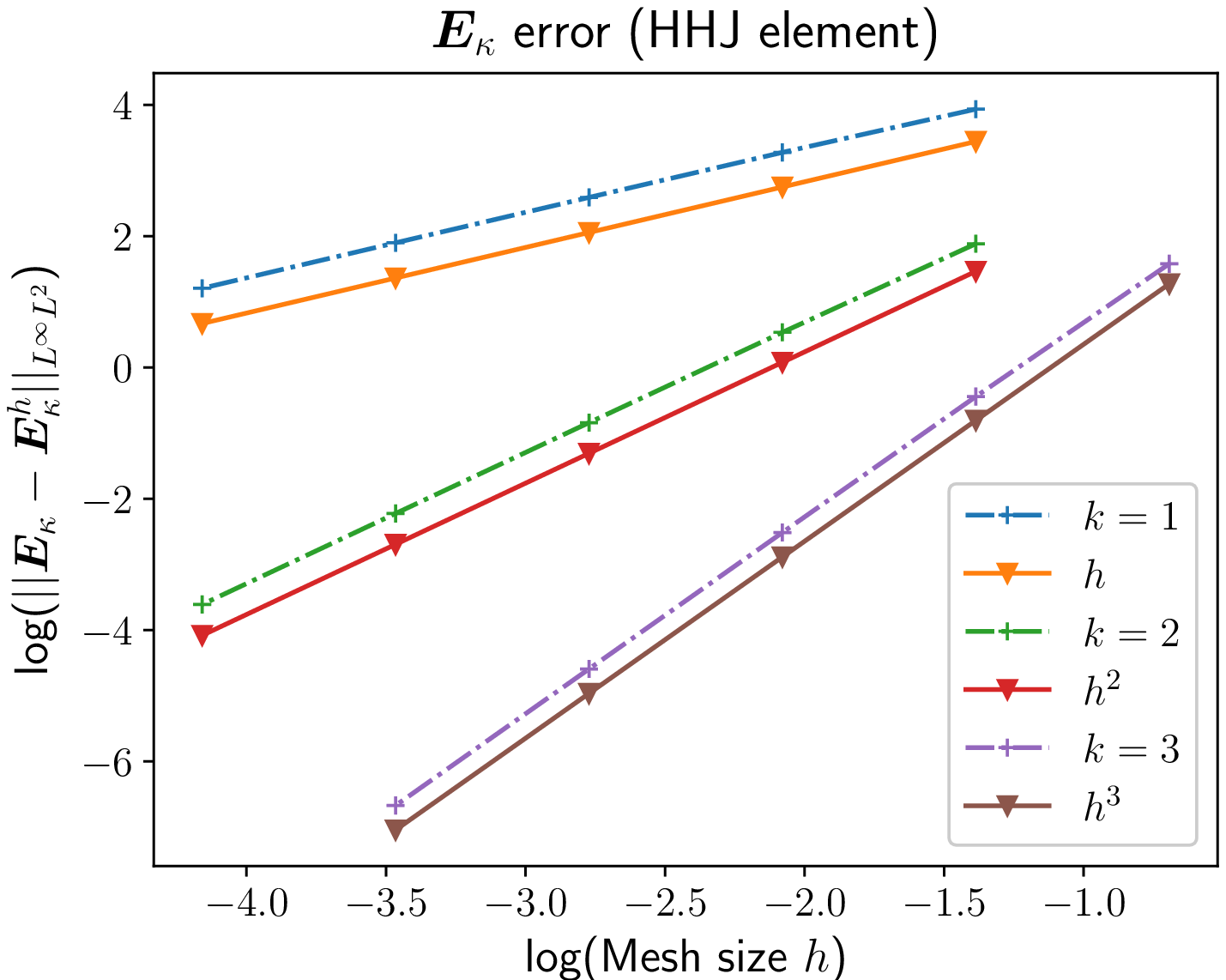}} \\
	\caption[errorHHF]{Error for the Kirchhoff plate using HHJ elements}%
	\label{fig:errorHHJ}%
\end{figure}

An analytical solution for the Kirchhoff plate is readily available. Consider the following solution of problem \eqref{eq:clKir} under simply supported conditions on a square unitary domain
\[
w^{\text{ex}}(x,y,t) = \sin(\pi x) \sin(\pi y) \sin(t), \quad  (x, y) \in (0,1)\times (0,1).
\] 
The forcing term is given by  
\[
f = (4 D \pi^4 - \rho b) \sin(\pi x) \sin(\pi y) \sin(t), \quad D = \frac{E_Y b^3}{12 (1-\nu^2)}.
\]
The corresponding variables in the port-Hamiltonian frame work are
\[
e_w^{\text{ex}} = \partial_t w^{\text{ex}}, \quad \bm{E}_\kappa^{\text{ex}} = \mathcal{D} \nabla^2 w^{\text{ex}}.
\]
Variables $(e_w^{\text{ex}}, \bm{E}_\kappa^{\text{ex}})$ under excitation $f$ solve problem~\eqref{eq:PH_sys_Kir_Ten}. The physical parameters used in simulation are reported in Table \ref{tab:parKir}. The weak form \eqref{eq:weak_kir_PH} and the finite elements \eqref{eq:HHJ} are considered. A direct solver with an LU preconditioner is used to compute the solution. Results are shown in Fig.~\ref{fig:errorHHJ}. The conjectured error estimates are respected.

\begin{table}[h]
	\centering
	\begin{tabular}{cccc}
		\hline 
		\multicolumn{4}{c}{Plate parameters} \\ 
		\hline 
		$E$ & $\rho$ & $\nu$  & $h$ \\
		136 $[\textrm{GPa}]$ & $5600\; [\textrm{kg}/\textrm{m}^3]$ & 0.3 &  0.001 $[\textrm{m}]$\\ 
		\hline 
	\end{tabular} 
	\captionsetup{width=0.95\linewidth}
	\vspace{1mm}
	\captionof{table}{Physical parameters for the Kirchhoff plate.}
	\label{tab:parKir}
\end{table}

\section{Conclusion}
In this paper, the link between mixed finite element method and pH plate models has been studied. It was shown that existing elements can be used to obtain structure-preserving discretization. A rigorous error analysis is still to be done but it should be easy to prove, given the available results. Since the pH framework provides a powerful description of boundary controlled systems, it is important that numerical methods be capable of handling generic boundary conditions. The methods discussed here possess this feature in the Mindlin plate case. For the Kirchhoff plate, a promising methodology is detailed in \cite{mixed_kirchhoff}, but the dynamical case has not been considered yet. Future developments include the analysis and discretization of viscoelastic and thermoelastic problems in pH form.

\begin{ack}
	The authors would like to thank Michel Sala\"un, Xavier Vasseur and Ghislain Haine from ISAE for the insightful and fruitful discussions.
\end{ack}

\bibliography{biblio_MTNS}             
                                                   







\appendix

\end{document}